\documentclass{article} 

\usepackage{fullpage} 
                     
\usepackage[frenchb,english]{babel}

\usepackage[T1]{fontenc} 
                        \usepackage{amsfonts}
\usepackage{amsmath} 
                     
\usepackage{graphicx}

\newtheorem{de}{Definition}[subsection] 
\newtheorem{theo}{Theorem}[section]    
\newtheorem{prop}[theo]{Proposition}    
             
\everymath{\displaystyle\everymath{}}

\author{Bertrand BARRAU }

\title{ Hilbert-Polya conjecture and demonstration of Riemann Hypothesis  
    \\
       \small }

\begin{document}
\maketitle

\begin{abstract}
 Using as starting point a classical integral representation of a L-function $L(s)$ we define a familly of two variables extended functions which are eigenfunctions of a Hermitian operator (having $i (s - \textstyle\frac{1}{2})$ as eigenvalues):

 $\mathbf{H}_{21} = \frac{1}{2} (\frac{\partial^2}{\partial^2 z}+ z^2) - i (y \frac{\partial}{\partial y} + \frac{1}{2})$  

This Hermitian operator can take also other forms, more symetric.

 In the case of particular L-functions, like Zeta function or Dirichlet L-functions, the eigenfunctions defined for this operator have symmetry properties.

Moreover, for $s$ zero fo Zeta function (or Dirichlet L-function), the associated eigenfunction has a specific property (a part of eigenfunction is canceled).
Finding such an eigenfunction, square integrable due to this "cancellation effect", would lead to Riemann Hypothesis using Hilbert-Polya idea.

\end{abstract}

\section{Introduction}   

Hilbert and Polya had the intuition that a Hermitian operator was "hidden" behind the Zeta function non trivial zeros\footnote{Zeros located on the "critical strip" of the complex plane defined by $0<Re(z)<1$.}: their imaginary parts  corresponding to the Hermitian operator eigenvalues. This hypothesis is closely linked to the Riemann Hypothesis\footnote {See \cite{TIT86} \cite{HME74}\cite{PAT88} for a good overview of the Zeta function theory and description of Riemann Hypothesis. See in Appendix a short reminder of Hilbert-Polya idea.} as eigenvalues of Hermitian operator are all real. 

With latest development of computing (millions of Zeta-function zeros are known) and probabilist studies \cite{BP01}, interesting distribution analogies were shown between Riemann Zeta function zeros and eigenvalues of some quantum mechanic operators, the intuition increased, and many conjectures were made about the form that could take the "hidden" Hermitian operator. A good list of references on the subject can be found in \cite{BK99}. 

These results and analogies support the Hilbert-Polya idea and motivates the search for Hermitian operators linked to Zeta function and to quantum mechanic\footnote{A summary of some latest approaches to this problem with their references are given in \cite{HR03}}.

If such an operator exists it would be interesting to know two of its main characteristics: its dimension and what is the origin of the "quantification effect" explaining the zero distribution on the critical line.

The first natural idea is to propose a one dimensional operator which will, by using the Sturm-Liouville theory, prove both alignment of zeros and their distribution\footnote{in this case "Quantification" of Zeta function zeros would be the results of classical  Sturm-Liouville theory in one dimension.}. Unfortunately, up to now, tentatives to find a one dimensional operator having its eigenvalues "explaining" the non trivial Zeta function zeros location have not succeed.

If it exists, the Hilbert-Polya operator could be a multidimensional operator with quantification effect not linked with Sturm-Liouville theory like it is proposed in some articles (see for example \cite{SO98} or \cite{GS08}).

 Article \cite{GS08} constructs a Hermitian operator "toy model" using mainly analogies and Zeta function zeros asymptotic properties\footnote{It should be outlined that this toy model operator, proposed by G. Sierra and P.K Townsend, is of the same type as the ones proposed in this article.}.

See also the study of Laplace-Beltrami Operator in \cite{MO97} which is also a two dimensional operator linked to Zeta function providing a spectral theory approach study but do not providing a spectral interpretation of the non trivial zeros.

The aim of this article is to show how the L-functions (and in particular Zeta function) are naturally linked to a two dimensional Hermitian operator taking different forms mainly as a combination of the operator $\frac{d^2}{dx^2}+x^2$ (named the IHO: Inverted Harmonic Oscillator) and the operator $x \frac{d}{dx}+\frac{1}{2}$.

The link between IHO and Zeta function has already been outlined in many publications, as for example in \cite{BK95} and \cite{BK99} where authors insisted on its importance and again in \cite{GS08} (In the form of the Harmonic Oscillator). The link of the operator $x \frac{d}{dx}+\frac{1}{2}$ to Zeta function is quite natural and has been outlined in different publications.\footnote {M. Berry and J. Keating for example proposed this operator as the Hilbert-Polya operator $H=xp$  \cite{BK98}.}
Moreover as outlined in \cite{BK98} and as it will be shown in this article, they fundamentaly correspond to the "same" operator as they are linked through a canonical transform.

We note: $\alpha=e^{-i \frac{\pi}{4}}$, more precisions on notations used in this article are given in Appendix.

\section{Definition of the extended functions $F_s(z)$ } 

We remark that to study the zeros of a L-function $L(1-s)=\displaystyle\sum_{n=1}^{\infty} \dfrac{a_n}{n^{1-s}}$ on the critical strip (we suppose the function convergent on this strip), it is equivalent to study the zeros of the following product: 

$L(1-s) \displaystyle\int_0^\infty  e^{-\pi  x^2} x^{-s} dx = \displaystyle\sum_{n=1}^{\infty} \dfrac{a_n}{n^{1-s}}\displaystyle\int_0^\infty  e^{-\pi  x^2} x^{-s} dx =   \displaystyle\int_0^\infty \displaystyle\sum_{n=1}^{\infty} a_n e^{-\pi n^2 x^2} x^{-s} dx $

It is well known that when the $a_n$ are such that a Poisson-Summation formula (or combination) exists for the function $\displaystyle\sum_{n=1}^{\infty} a_n e^{-\pi n^2 x^2}$, this type of relation is used to show that $L(s)$ verifies a functional equation.

In the case of L-function not converging in the critical strip, the study of the zeros can be done by studying the above integral representation slightly modified. 

For example in the case of the Zeta function, as  $\displaystyle\sum_{n=1}^{\infty} \dfrac{1}{n^s}$ does not converge for $Re(s)<1$, the precedent integral representation can not be used as such to study the Zeta zeros.

The relation above has to be slightly modified to avoid problem of integral convergence in zero (See \cite{HME74} p213) and can take the form\footnote{$\xi(s)=(s-1) \pi^{-\frac{s}{2}} \Gamma(\frac{s}{2}) \zeta(s)$}:

\vspace {0.3 cm}

$\dfrac{2 (s-1) \pi^{-\frac{s}{2}} \Gamma(\frac{s}{2}) \zeta(s)}{s (s-1)}=\dfrac{2 \xi(s)}{s (s-1)}=\displaystyle\int_0^\infty {(\displaystyle\sum_{n \neq 0} e^{-\pi n^2 x^2}-\dfrac{1}{x}) x^{-s} dx}$       \; \;                      for $0 < Re(s) < 1$ \; \; (1)

And in this case the relation $G(x)=\displaystyle\sum_{n=-\infty}^{\infty} e^{-\pi n^2 x^2} = \dfrac{1}{x} \displaystyle\sum_{n=-\infty}^{\infty} e^{-\pi \frac{n^2}{ x^2}}$ is used to show the functional equation of $\zeta(s)$.

 Using the change of variable $x \to \frac{1}{x}$ we find the functional equation:

$\dfrac{2 \xi(s)}{s (s-1)}= \displaystyle\int_0^\infty {(\displaystyle\sum_{n=-\infty}^{\infty} e^{-\pi \frac{n^2}{ x^2}}-1-x) x^{s-1} \frac{dx}{x}}=\displaystyle\int_0^\infty {(\displaystyle\sum_{n \neq 0} e^{-\pi n^2 x^2}-\dfrac{1}{x}) x^{s-1} dx}=\dfrac{2 \xi(1-s)}{s (s-1)}$ \; \; (2)

And subtracting relation (1) to the same one obtained following a change of variable $ x \to 2x$ and multiplication by $2^s$, we obtain:

$\dfrac{2 \xi(s)}{s (s-1)}(1-2^s)=\displaystyle\int_0^\infty {\displaystyle\sum_{n \neq 0}  (-1)^{n-1}e^{-\pi n^2 x^2} x^{-s} dx}$    \; \;  

So we see that for the study of Zeta function zeros we can use the $\eta(s)$ L-function\footnote{$\eta(s)= \displaystyle\sum_{n=1}^{\infty} \dfrac{(-1)^{n-1}}{n^s}$ \; \; this function as same zeros as Zeta on the critical strip.} to be in the case of a converging L-function on the critical strip as presented at the beginning of this paragraph.

It has to be noted that integral representations as above can be done with other functions than $e^{-n^2 x^2 }$, but the use of this specific function (having good convergence and Fourier properties) in the integral representation of $L(s)$ is particularly relevant as when the $a_n$ have specific properties (for example when they are a Dirichlet Character, or all identical to 1 as in Zeta), the functional equation of $L(s)$ (using Poisson summation formula and change of variable  $x \to \frac{1}{x}$) can easily be demonstrated as already shown above (2).

The object of this article is to show how these type of integral representation of L-function can be naturally extended to functions of one and then two variables which are eigenfunctions of a Hermitian operator which can take different forms.

The first step will be to extend the integral representation by using the function  $e^{\pi (z-i n x)^2 }$ instead of simply $e^{-n^2 x^2 }$, we will see that this extension which is natural leads to a one dimensional operator which makes the imaginary part of $s$ appearing as eigenvalue.

This new function family of $z$, we will be noted $H_s^{L}(z)$ as we will see they are solution of the "Hermite polynomial" differential equation.

We define $a_{-n}=a_n$, this artificial extension is useful as Poisson formula (which is often used) necessitates a sum from $-\infty$ to $+\infty$.

\begin{de}
 
 For $L(s)=\displaystyle\sum_{n=1}^{\infty} \dfrac{a_n}{n^s}$ \; we define:  $H^{L}_s(z)=\displaystyle\int_0^\infty \displaystyle\sum_{n \neq 0} a_n e^{\pi (z-i n x)^2 }  x^{-s} dx $

And:

 for Zeta function:    $H^\zeta_s(z)=\displaystyle\int_0^\infty (\displaystyle\sum_{n \neq 0} e^{\pi (z-i n x)^2 }-\dfrac{1}{x}) x^{-s} dx $

\end{de}

So for Eta function we have: $H^\eta_s(z)=\displaystyle\int_0^\infty \displaystyle\sum_{n \neq 0} (-1)^{n-1} e^{\pi (z-i n x)^2 } x^{-s} dx $

We see $H^{L}_s(z)$ is an extension of the initial relation and that in particular for $z=0$ we come back to the original not extended integral relation: 

$H^{L}_s(0)=\displaystyle\int_0^\infty \displaystyle\sum_{n \neq 0} a_n e^{-\pi n^2 x^2 }  x^{-s} dx =\pi^{\frac{s-1}{2}} \Gamma(\frac{1-s}{2})  L(1-s)$

And in case of Zeta function:
 $H^\zeta_s (0)=\displaystyle\int_0^\infty {(\displaystyle\sum_{n \neq 0} e^{-\pi n^2 x^2}-\dfrac{1}{x}) x^{-s} dx}=\dfrac{2 \xi(s)}{s (s-1)}$

\vspace {0.3 cm}

This definition is valid for all real $z$ and for $s$ in the critical strip (as we assume the sum defining $L(s)$ converges on this strip).

We will now see that this extension has the specific property to be solution of a differential equation and to be identically nul for $s$ zero of Zeta-function.

\section{ Extended $H^{L}_s(z)$ functions are eigenfunction of the Hermite differential operator} 

\begin{prop}\label{Theorem 1: The symmetrical function differential equation} 

Considering the function $H^{L}_s(z)=\int_0^\infty \sum_{n \neq 0}  a_n e^{\pi (z-i n x)^2 } x^{-s} dx$ \; 

and  $H^{\zeta}_s(z)=\int_0^\infty (\sum_{n \neq 0}  e^{\pi (z-i n x)^2 }-\dfrac{1}{x}) x^{-s} dx$ \; 

we have for $ 0 < Re(s) < 1$ : 

 $H^{\zeta}_s(z)''- 2\pi zH^{\zeta}_s(z)'= 2\pi s H^{\zeta}_s(z)$ and $H^{L}_s(z)''- 2\pi zH^{L}_s(z)'= 2\pi s H^{L}_s(z)$

 \end{prop}

Verification of this property is obvious for $H^{L}_s(z)$ as we can "extract" the function $L(s)$ from the integral to write:

$H^{L}_s(z)= 2 L(1-s) \int_0^\infty  e^{\pi (z-i  x)^2 } x^{-s} dx$ \; 

And then easily check that this function of $z$ verifies the differential equation.

In the case of $H^{\zeta}_s(z)$ the Zeta-function can not be extracted and a direct calculation provides the result. (Splitting the integral in two parts and using Poisson summation formula to avoid any problem of integral convergence in zero)

Note that the function $(\sum_{n \neq 0}  e^{\pi (z-i n x)^2 }-\dfrac{1}{x})$ under the integrand has a finite limit for $x$ tending to zero  as we have following Poisson summation formula for $z$ fixed and $x$ real not null:

$\sum_{n=-\infty}^{\infty}e^{- \pi x^2  n^2 - 2 i \pi n  z x +\pi z^2}= \frac{1}{x} \sum_{n=-\infty}^{\infty}e^{-\pi ( \frac{n}{x}-z)^2+\pi z^2 }$ and so:

$\lim_{x \to 0} (\sum_{n \neq 0}  e^{\pi (z-i n x)^2 }-\dfrac{1}{x})=\lim_{x \to 0}  \frac{1}{x} \sum_{n \neq 0} e^{-\pi ( \frac{n}{x}-z)^2+\pi z^2 }-e^{\pi z^2} = -e^{\pi z^2}$ 

The good property of this function and the fact that $\sum_{n \neq 0}  e^{\pi (z-i n x)^2 }$ is going exponentially to 0 at infinity will ensure convergence of integrals and correctness of our future manipulations.

\vspace {0.3 cm}

We see that $2 \pi s$ is eigenvalue of eigenfunction $H^{L}_s(z)$, at this point we make a natural change of variable to this differential equation (the one done to pass from Hermite polynomials to Hermite functions), and we see that then changing $x$ by $x e^{-i \frac{\pi}{4}}$ makes $\lambda$ appearing as eigenvalue ($s=\frac{1}{2} + i \lambda$).

The term $\frac{1}{2}$ required for the apparition of $\lambda$ appears naturally.

\begin{prop} We define the  $\varphi^{L}_s(z)$ functions by:

$\varphi^{L}_s(z)=e^{\frac{i}{2}z^2} H^{L}_s(z)(\frac{z}{\sqrt{\pi}} e^{-i \frac{\pi}{4}})$   

  and we have the following properties (Posing $\alpha=e^{-i \frac{\pi}{4}}$):

\begin{itemize}

\item $\varphi^{L}_s(z)= e^{ i \frac{z^2}{2}} \int_0^\infty  \sum_{n \neq 0}^{} a_n e^{ (\alpha z-i  \sqrt{\pi} n x)^2 }  x ^{-s}dx$ ; 

\item $\varphi^{\zeta}_s(z)= e^{ i \frac{z^2}{2}} \int_0^\infty (\sum_{n \neq 0}^{}e^{ (\alpha z-i  \sqrt{\pi} n x)^2 }-\dfrac{1}{x}) x ^{-s}dx$ ;

\item $\varphi^{\zeta}_s(z)''+  z^2  \varphi^{\zeta}_s(z) = i(1-2 s)  \varphi^{\zeta}_s(z) = 2 \lambda\varphi^{\zeta}_s(z)$

\item $\varphi^{L}_s(z)''+  z^2  \varphi^{L}_s(z) = i(1-2 s)  \varphi^{L}_s(z) = 2 \lambda\varphi^{L}_s(z)$

\end{itemize}
\end{prop}

We obtain these results for $\varphi^{L}_s(z)$ by simple calculation following the change of variable proposed and properties of $H^{L}_s(z)$ already shown.

We see how the $H_s$ is a sort of analog to the $H_n$ (Hermite polynomials) and the $\varphi_s$ an analog to the $\varphi_n$ (Hermite functions)\footnote{Link between Hermite functions and Zeta function as already been noticed in different other ways see for example \cite{DB97}, this analogy can be developped: generating function, Fourier properties...}

We can easily check that (as $L(s)$ can be extracted from integral),  $\varphi^{L}_s(z)=0$ for all z if and only if $s$ is zero of $L(1-s)$ (So also zero of  $L(s)$ when $a_n$ are such that functional equation provides a relation between symmetry $L(1-s)$ and $L(s)$).

Identically we can prove that $\varphi^{\zeta}_s(z)=0$ for all z if and only if $s$ is zero of $\zeta(s)$:

 $\varphi^{\zeta}_s(z)$  is in the one dimensional vector space of even solutions of differential equation above. 
On other hand, the even solutions of the differential equation are the even Parabolic Cylinder functions, and with notation of \cite{MAIS72} for $y_0(z)$ the even Parabolic Cylinder functions solution of IHO\footnote{IHO=Inverted Harmonic Oscillator} such that  $y_0(0)=1$ (Taking notations of \cite{MAIS72} Page 692), we have, (for $a=\lambda$ and as $\varphi^{\zeta}_s(0)= \dfrac{2 \xi(s)}{s (s-1)}$):

$\varphi^{\zeta}_s(z)= \dfrac{2 \xi(s)}{s (s-1)} y_{0}(\sqrt{2}z)$ 

At this stage we can wonder about the interest of this extension as the L-function can still be extracted from its integral and $H^{L}_s(z)$ written as a pure product involving $L(1-s)$ and an integral.

But the specificty of this extension is that it respects the functional equation, using Poisson Summation formula:

$\sum_{n=-\infty}^{\infty}e^{- \pi x^2  n^2 - 2 i \pi n \alpha z x }= \frac{1}{x} \sum_{n=-\infty}^{\infty}e^{-\pi ( \frac{n}{x}-\alpha z)^2 }$ and change of variable $x \to \frac{1}{x}$ we have:

$\varphi^{\zeta}_s(z)=e^{ i \frac{z^2}{2}} \int_0^\infty (\sum_{n \neq 0}^{}e^{ (\alpha z-i  \sqrt{\pi} n x)^2 }-\dfrac{1}{x}) x ^{-s}dx = e^{- i \frac{z^2}{2}} \int_0^\infty (\sum_{n \neq 0}^{}e^{ (i \alpha z-i  \sqrt{\pi} n x)^2 }-\dfrac{1}{x}) x ^{-s}dx= \varphi^{\zeta}_{1-s} (iz) $

So finally:  $\varphi^{\zeta}_s(z)= \varphi^{\zeta}_{1-s} (iz) $

Identical relations can be found for $\varphi^{\eta}_s(z)$ and for $\varphi^{L}_s(z)$ providing $a_n$ have sufficient symmetry properties: for example in the case $a_n$ is a character $\chi_n$.

We can also write $\varphi^{\zeta}_s(z)$ in the following way to see its interesting symmetrical properties due to Poisson summation formula (Splitting integral in two and changing $x$ by$\textstyle\frac{1}{x}$ in first integral):

$\varphi^{\zeta}_s(z)=\displaystyle\int_1^\infty (\sum_{n \neq 0}^{}e^{ -i \frac{z^2}{2} - 2i\sqrt{\pi}\alpha nx z-\pi n^2 x^2 }-\dfrac{e^{ i \frac{z^2}{2}}}{x}) x^{-s}dx +\int_1^\infty (\displaystyle\sum_{n \neq 0}^{}e^{ i \frac{z^2}{2} - 2\sqrt{\pi}\alpha nx z-\pi n^2 x^2 }-\dfrac{e^{- i \frac{z^2}{2}}}{x}) x^{s-1} dx$

From this expression we deduce that $\varphi^{\zeta}_s(z)$ takes real values for $s$ on the critical line (so the integral multiplying the Zeta function is particularly adapted as it "twisted" the Zeta function on the critical line to make it real).

Remark: When considering Zeta and Eta functions, we see that the associated functions we defined can also be defined using Theta functions:

with $\theta_3(\tau,y)=\sum_{n=-\infty}^{\infty} e^{ i \pi \tau n^2 + 2i\pi n y}$  (notation of \cite{YH97}), the expression of $\varphi_s^{\zeta}(z)$ becomes:

\vspace {0.3 cm}
$\varphi_s^{\zeta}(z)=e^{i \frac{z^2}{2}} \int_0^\infty ( \theta_3(i x^2, -x z \dfrac{\alpha}{\sqrt{\pi}})-1-\dfrac{1}{x}) x^{-s}dx= e^{-i \frac{z^2}{2}} \int_0^\infty ( \theta_3(i x^2, -i x z \dfrac{\alpha}{\sqrt{\pi}})    -1-\dfrac{1}{x}) x^{s-1}dx$

\vspace {0,3 cm}

In the same way we can use  $\theta_2$ and  $\theta_4$ to express eigenfunctions of $\eta$ function:

$\varphi_s^{\eta}(z)=e^{i \frac{z^2}{2}} \int_0^\infty ( \theta_2(i x^2, -x z \dfrac{\alpha}{\sqrt{\pi}}) -1) x^{-s}dx= e^{-i \frac{z^2}{2}} \int_0^\infty ( \theta_4(i x^2, -i x z \dfrac{\alpha}{\sqrt{\pi}})    -1) x^{s-1}dx$

(Using fact that $ \theta_2(x,z)= \frac{1}{x} \theta_4( x,i z) $ and change of variable $x \to \frac{1}{x}$ as usual in this article.)

\section {Two dimensional Hermitian operators associated to L-functions } 

In precedent paragraph we have seen that $\varphi^{L}_s(z)$ functions are eigenfunctions of the Hermitian operator which is the IHO. We have also seen that in this form there is still no "actual result" as the  $\varphi^{L}_s(z)$ functions are in fact proportional to the $L(1-s)$ functions (the $L(s)$ functions "escape" from the integral with their intrinsic properties). In order to "keep" the L-function properties inside the integrand we need to add another variable. This should be done keeping the Hermitian property of the operator acting on the new eigenfunction and keeping $\lambda$ as eigenvalue. 

The interesting result is that such extension under these two constraints is possible.

We start with eigenfunctions of two dimensional operators (which are extension of $\varphi^{\zeta}_s(z)$) having $s$ as eigenvalues.

\begin{de} 
For any fixed "good" function $g(t) \in \mathcal{C}^2 (\mathbb{R^+},\mathbb{C})$, 
we define a family of functions $A^{L}_{s,a,g}(z,y)$ on $\mathbb{R}\times \mathbb{R^+}$   for $a \neq 0$ by:

 $A^{L}_{s,a,g}(z,y)= \int_0^\infty (\sum_{n \neq 0}^{} a_n e^{\pi ( z-i n x)^2 }) g(xy^a) x^{-s}dx$

 $A^{\zeta}_{s,a,g}(z,y)= \int_0^\infty (\sum_{n \neq 0}^{}e^{\pi ( z-i n x)^2 }-\dfrac{1}{x}) g(xy^a) x^{-s}dx$

\end{de}

For functions $g(t)$ "good enough" we obtain functions $A_{s,a,g}$ well defined and eigenfunction of an operator depending only on $a$:

\begin{prop}\label{Theorem: The main equation} 
:

On the domain they are well defined and in $\mathcal{C}^2(\mathbb{R}) \times \mathcal{C}^1 (\mathbb{R}^+)$, the $A_{s,a,g}(z,y)$  are eigenfunctions of the fundamental  $P_a$ operator:

 $P_{a}= \frac{1}{ 2\pi} {\dfrac{\partial^2}{\partial^2z}} -  z \dfrac{\partial }{\partial z} + \frac{y}{a} {\dfrac{\partial }{\partial y}} $ 

And we have:  $P_{a} A_{s,a,g}(z,y)=   s A_{s,a,g}(z,y)$

\end{prop}
Demonstration is a direct calculation.

This is an intermediate result: we need the same change of variable as above to have a two dimensional Hermitian operator having $\lambda$ as eigenvalue. 

\begin{prop} 

For $s$ in the critical strip, we define the $B_{s}(z,y)$ functions from $ \mathbb{R} \times \mathbb{R^*}^+  $ to $\mathbb{C}$ by\footnote{We remind that $\alpha=e^{-i \frac{\pi}{4}}$ and $s = \scriptstyle\frac{1}{2}+i \lambda$}:
\vspace {0.5 cm}

 $B^{L}_{s,a}(z,y)= \frac{e^{\frac{i}{2}z^2}}{\sqrt{y}} A^{L}_{s,a,g}(\frac{z}{\sqrt{\pi}} e^{-i \frac{\pi}{4}}, y) $ ; 

So: 

$B^{L}_{s,a}(z,y)=e^{ i \frac{z^2}{2}}y^{-\frac{1}{2}} \int_0^\infty \sum_{n \neq 0}^{} a_n e^{ (\alpha z-i  \sqrt{\pi} n x)^2 } g(xy^a) \; x ^{-s}dx  $

and: $B_{s,a}^{\zeta}(z,y)= \frac{e^{ i \frac{z^2}{2}}}{ \sqrt{y}} \int_0^\infty (\sum_{n \neq 0}^{}e^{ (\alpha z-i  \sqrt{\pi} n x)^2 }-\dfrac{1}{x}) g(xy^a) \; x ^{-s}dx$

 Noting $\mathbf{H}_{21}^{a}=\dfrac{1}{2} ({\dfrac{\partial^2}{\partial^2z}} + z^2 )- \frac{i}{a} (y{\dfrac{\partial }{\partial y}}+\dfrac{1}{2}) $ we have on $ \mathbb{R} \times \mathbb{R^*}^+  $:

 \begin{itemize}

\item  $\mathbf{H}_{21}^{a} B^{L}_{s,a}(z,y)=\lambda   B^{L}_{s,a}(z,y)$;

\end{itemize}
\end{prop}

This is the immediate result of a variable change in the differential equation verified by $A_s$ functions given above.

The notation has been choosen because the differential operator $\mathbf{H}_{21}^{a}$ is of second order in $z$ and first order in $y$.

 It is interesting to note that this operator does not depend on the L-function.

 The operator  $\mathbf{H}_{21}^{a}=\dfrac{1}{2} ({\dfrac{\partial^2}{\partial^2z}} + z^2 )-\frac{ i}{a} (y{\dfrac{\partial }{\partial y}}+\dfrac{1}{2}) $ is hermitian on $\mathcal{H}_0$ (Classical result: by integration by parts we obtain $<\mathbf{H}_{21} f,  g> = <f, \mathbf{H}_{21} g>$).

\section {On the different forms taken by the Hermitian operator associated to L-functions} 

From now we will take $a=1$ in order to ease the notation (but all what follows can be consider with any value of $a \neq 0$ without difficulty) 

We will now define the transform\footnote{which is a variant of Weierstrass transform.} allowing to pass from operator $i (x \frac{d}{dx} + \frac{1}{2})$ to operator $ \frac{1}{2} ( \frac{d^2}{dx^2} +x^2 ) $

\begin{prop} Considering  $f(t)$ a function such that: $\lim_{t \to -\infty} tf(t)e^{ - t^2 }=0$ and  $\lim_{t \to \infty} tf(t)e^{  - t^2 }=0$

and considering the transform $T$: $T(f)(x)= e^{   - \frac{x^2}{2} }  \int_{-\infty}^\infty e^{ ( x-it)^2} f(t) dt$

If $  \phi (x) = i (x \frac{d}{dx}  + \frac{1}{2} ) f(x)$

Then $T(\phi )(x) =  - \frac{i}{2} ( \frac{d^2}{dx^2}  - x^2 ) T(f)( x) $ and $T(\phi )(\alpha x) =  \frac{1}{2} ( \frac{d^2}{dx^2}  +x^2 ) T(f)(\alpha x) $

\end{prop}

This result is a simple calculation with integration by parts.

We see here that this transform appears in the expression of the eigenfunctions $A_s$ and $B_s$. 
And that using this transform the hermitian operator $\mathbf{H}$ we defined and its eigenfunctions can take following forms:

$\mathbf{H}_{11}= i(x \frac{\partial}{\partial x} - y \frac{\partial}{\partial y}) \longleftrightarrow    x^{-s} y^{-\frac{1}{2}}  \sum_{n=1}^{\infty} \dfrac{a_n}{n^{1-s}} g( \frac{  x y}{ n})   $

$\downarrow T_x$

$\mathbf{H}_{21}= \frac{1}{2} (\frac{\partial^2}{\partial^2 z}+ z^2) - i (y \frac{\partial}{\partial y} + \frac{1}{2}) \longleftrightarrow  B_s(z,y)=e^{ i \frac{z^2}{2}}y^{-\frac{1}{2}} \int_0^\infty \sum_{n \neq 0}^{} a_n e^{ (\alpha z-i  \sqrt{\pi} n x)^2 } g(xy) \; x ^{-s}dx  $

$\downarrow T_y$

$\mathbf{H}_{22}= \frac{1}{2}( \frac{\partial^2}{\partial^2 z}+ z^2) -\frac{1}{2}( \frac{\partial^2}{\partial^2 u}+ u^2)  \longleftrightarrow  C_s(z,u)= e^{ i  \frac{z^2}{2}} e^{ i  \frac{u^2}{2}} \int_0^\infty \int_0^\infty \sum_{n \neq 0}^{} a_n  e^{  (\alpha z-i \sqrt{\pi} n x)^2 } e^{  (\alpha u-i \sqrt{\pi}  y)^2 } g(xy) y^{-\frac{1}{2}}x ^{-s} dx dy$   

These three forms of the same operator with their eigenfunctions have all  eigenvalues $\lambda$ (the imaginary part of $s$) linked to the associated L-functions.
In the three cases the operator associated is Hermitian, and to be able to use the idea of Hilbert-Polya, it is required to show that there exists one function $g(t)$ such that when $s$ is zero of $L(s)$ one of the above eigenfunction is in the associated Hilbert space. (i.e. eigenfunction is module square integrable and null at edge of integration domain considered for Hilbert product).
This condition is obviously impossible for the operator in its "initial form" $\mathbf{H}_{11}$, this operator will never have square module integrable eigenfunctions. (In the same way Hermitian operator $i(x \frac{\partial}{\partial x} + \frac{1}{2})$ as no eigenfunction in $L^2(\mathbb{R}^+$)

\begin{prop}
If there exists a function $g(t)$ such that for $s$ zero of Zeta function we have\footnote{We note $\mathcal{H}_0 =\{ h(z,y) \in \mathbb{C};  h \in L^2 (\mathbb{R} \times \mathbb{R}^+) ; \forall z, h(z,0)=0\}$ and $\mathcal{H}_1 =\{ h(z,y) \in \mathbb{C};  h \in L^2 (\mathbb{R} \times \mathbb{R}) \}$}:

$e^{ i \frac{z^2}{2}}y^{-\frac{1}{2}} \int_0^\infty \sum_{n \neq 0}^{} (-1)^{n-1}  e^{ (\alpha z-i  \sqrt{\pi} n x)^2 } g(xy) \; x ^{-s}dx \in \mathcal{H}_0$

or $e^{ i  \frac{z^2}{2}} e^{ i  \frac{u^2}{2}}  \int_0^\infty \int_0^\infty \sum_{n \neq 0}^{}  (-1)^{n-1}   e^{  (\alpha z-i \sqrt{\pi} n x)^2 } e^{  (\alpha u-i \sqrt{\pi}  y)^2 } g(xy) y^{-\frac{1}{2}}x ^{-s} dx   \in \mathcal{H}_1$ 

then all the zeros of Zeta functions have their real part equal to $\frac{1}{2}$.

\end{prop}

Results is immediate using the fact that Operator is Hermitian, as we have already seen, if $B_s$ is in $\mathcal{H}_0$ for example:

$\lambda <B_{s} , B_{s}>=<\mathbf{H}_{21} (B_{s}) , B_{s}> = \int_{- \infty}^{\infty} \int_{0}^{\infty} ( \dfrac{1}{2} ({\dfrac{\partial^2}{\partial^2z}}  B_{s}+ z^2  B_{s})-i(y{\dfrac{\partial }{\partial y}}B_{s}+\dfrac{1}{2}  B_{s})) \; \overline{ B_s(z,y)}  dy dz= <B_{s} ,\mathbf{H}_{21} (B_{s})>= \overline{\lambda} <B_{s} , B_{s}>$

Now we will explain why there is a chance that such a function $g(t)$, which makes $B^{\eta}_s$ module square integrable for $s$ zero of Zeta, exists:

Taking the function $B^{\eta}_s$ we know that when $s$ is zero of zeta function then for all $z$:

 $e^{ i \frac{z^2}{2}}y^{-\frac{1}{2}} \int_0^\infty \sum_{n \neq 0}^{} (-1)^{n-1} e^{ (\alpha z-i  \sqrt{\pi} n x)^2 }  \; x ^{-s}dx  =0$

Meaning that in this case ($s$ zero of zeta) and only in this case we have whatever the constant A is:

$B_s(z,y)=e^{ i \frac{z^2}{2}}y^{-\frac{1}{2}} \int_0^\infty \sum_{n \neq 0}^{} (-1)^{n-1} e^{ (\alpha z-i  \sqrt{\pi} n x)^2 }  (g(xy)-A)\; x ^{-s}dx  = e^{ i \frac{z^2}{2}}y^{-\frac{1}{2}} \int_0^\infty \sum_{n \neq 0}^{} (-1)^{n-1} e^{ (\alpha z-i  \sqrt{\pi} n x)^2 }  g(xy) \; x ^{-s}dx$

\textbf{For example} taking $g(t)= e^{-\frac{1}{t^2}}$ we can write $B_s$ (eigenfunction of $\mathbf{H}_{21}$) in two different ways:

$\begin{cases} 
B_s(z,y)=e^{ i \frac{z^2}{2}}y^{-\frac{1}{2}} \displaystyle\int_0^\infty \displaystyle\sum_{n \neq 0}^{} (-1)^{n-1} e^{ (\alpha z-i  \sqrt{\pi} n x)^2 }  e^{-\frac{1}{x^2y^2}} \;  x ^{-s}dx  \; \; \; \text{  (a)}
\\
B_s(z,y)=e^{ i \frac{z^2}{2}}y^{-\frac{1}{2}} \displaystyle\int_0^\infty \displaystyle\sum_{n \neq 0}^{} (-1)^{n-1} e^{ (\alpha z-i  \sqrt{\pi} n x)^2 }  ( e^{-\frac{1}{x^2y^2}}-1) \; x ^{-s}dx  \; \; \;  \text{(b)}
\end{cases}$

To study if $B_s(z,y)$ is module square integrable, we need to assess its behavior in $+\infty$ and $0$ in $z$ and $y$. 

For $y$ tending to 0 we use expression (a) which shows that $B_s(z,y)$ should tend to zero (as in this case $e^{-\frac{1}{x^2y^2}}$ tends to zero) and for $y$ tending to infinity we use the expression (b) which should tend to zero (as in this case $( e^{-\frac{1}{x^2y^2}}-1)$ tends to zero). 

The condition $s$ zero of Zeta function cancels the constant term of $g(t)$, and then there is a possibility that this cancellation involves the square integrable condition (in the same way $e^{-\frac{1}{t^2}}$ is not in $L^2(\mathbb{R}^+)$ but the function defined near zero by $e^{-\frac{1}{t^2}}$ and near $+\infty$ by $e^{-\frac{1}{t^2}}-1$ is in $L^2(\mathbb{R}^+)$). (This cancellation property can be found directly by using Laplace method to determine asymptotic limit of $B_s$ for $z$ tending to infinity: we see that the first term of development disappears if an only if $s$  is zero of Zeta function)

Note that $B_s$ takes also following expression (by change of variable and Poisson formula\footnote {$\sum_{n \in \mathbb{Z}} (-1)^{n-1} e^{i \frac{z^2}{2}+( \alpha z - i \sqrt{\pi} n x )^2} =  \frac{1}{x} \sum_{n \in \mathbb{Z}} e^{ -i \frac{z^2}{2} -(\frac{2n+1}{2x}\sqrt{\pi} - \alpha z)^2} $}):

$B_s(z,y)= e^{ i \frac{z^2}{2}}y^{-\frac{1}{2}} \int_0^\infty ( \sum_{n \in \mathbb{Z}}^{} (-1)^{n-1} e^{ -(\frac{2n+1}{2} x \sqrt{\pi} - \alpha z)^2}-1) e^{- \frac{ x^2}{ y^2}} x^{s-1}dx$

\vspace {0.3 cm}

\textbf{Why for some L-function there is no chance to find eigenfunctions module square integrable:}

\vspace {0.1 cm}

One more condition to find a good function $g$ is to have the $a_n$ specific, as already explained, to be able to use the change of variable $x \to \frac{1}{x}$ and Poisson Summation formula as already mentioned: this manipulation allows to change variable of $g$ from $x$ to $\frac{1}{x}$ while the rest of the integral remains nearly unchanged.
 
The L-functions with associated functional equations of the form: $\Lambda (s)= r_s  \Lambda (k-s)$ with $k>1$ will not have $a_n$ coefficient compatible with Poisson Summation formula involving $e^{\pi(z-ix)^2}$, therefore these L-function (from automorphic forms for example) even if having eigenfunction for above operators, will never have systematically their eigenfunctions associated to their zero in $L^2$, and this is compatible with the fact that L-functions associated to automorphic forms do not have all their zeros on the critical line.

\vspace {0.3 cm}
\textbf{Eigenfunctions candidates to be module square integrable}
\vspace {0.1 cm}

For zeros of the Zeta-Function (as for Eta-function), the immediate eigenfunctions candidates that could be in $L^2$ when $s$ is zero (so when we have the specificity explained above: the constant term disappearing under the integrand) could be: 

\textbf{For Operator $\mathbf{H}_{21}$:}

Taking $g(t) =  e^{-\frac{k}{t^2}} $  for $k$ complex: $B_s(z,y)=e^{ i \frac{z^2}{2}}y^{-\frac{1}{2}} \displaystyle\int_0^\infty \displaystyle\sum_{n \neq 0}^{} (-1)^{n-1} e^{ (\alpha z-i  \sqrt{\pi} n x)^2 }  e^{-\frac{k}{x^2y^2}} \;  x ^{-s}dx $

Another interesting idea is to consider eigenfunction obtain by taking $g(t)= \sum_{p \in \mathbb{Z}} (-1)^{p} e^{ - \pi p^2 t^2}$

then, using the:   $\sum_{n \in \mathbb{Z}} (-1)^{n} e^{ - \pi n^2 x^2 y^2} = - \frac{1}{xy} \sum_{n \in \mathbb{Z}} e^{- \pi  \frac{(2n+1)^2}{4x^2 y^2}} $

and using that in this case the cancellation of constant term explained above is done by first term on second term but also by second term of integrand on first one as we have:

$ \int_0^\infty  (\sum_{p \in \mathbb{Z}} (-1)^{p} e^{ - \pi p^2 x^2 y^2} -1)   x^{-s}dx = 0 $

using these properties we see that the eigenfunction takes the following forms if and only if $s$ is zero of Zeta (Initial expression of $B_s$ in this case is given by (3)):

\vspace {0.2 cm}
$ B_s = e^{ i \frac{z^2}{2}}y^{-\frac{1}{2}}  \int_0^\infty  (\sum_{n \neq 0}  (-1)^{n} e^{ (\alpha z-i \sqrt{\pi}n \frac{1}{x})^2 }  )( \frac{1}{y} \sum_{p \in \mathbb{Z}} e^{- \pi  \frac{(2p+1)^2 x^2}{4 y^2}} )     x^{s-1}dx$

$  = e^{ i \frac{z^2}{2}}y^{-\frac{1}{2}}  \int_0^\infty  (\sum_{n \neq 0}  (-1)^{n} e^{ (\alpha z-i \sqrt{\pi}n x)^2 }  )( \frac{1}{xy} \sum_{p \in \mathbb{Z}} e^{- \pi  \frac{(2p+1)^2}{4x^2 y^2}} )     x^{-s}dx$

$  = e^{ i \frac{z^2}{2}}y^{-\frac{1}{2}}  \int_0^\infty  (\sum_{n \neq 0}  (-1)^{n} e^{ (\alpha z-i \sqrt{\pi}n x)^2 }  ) ( \sum_{p \in \mathbb{Z}} (-1)^{p} e^{ - \pi p^2 x^2 y^2})   x^{-s}dx$ \; \; (3)

$  =e^{ i \frac{z^2}{2}}y^{-\frac{1}{2}}  \int_0^\infty  (\sum_{n \neq 0}  (-1)^{n} e^{ (\alpha z-i \sqrt{\pi}n x)^2 }  ) (\sum_{p \in \mathbb{Z}} (-1)^{p} e^{ - \pi p^2 x^2 y^2} -1)   x^{-s}dx$

$  =e^{ i \frac{z^2}{2}}y^{-\frac{1}{2}}  \int_0^\infty  (\sum_{n \neq 0}  (-1)^{n} e^{ (\alpha z-i \sqrt{\pi}n x)^2 }  + e^{-iz^2} ) (\sum_{p \in \mathbb{Z}} (-1)^{p} e^{ - \pi p^2 x^2 y^2} -1)   x^{-s}dx$

$  =e^{ i \frac{z^2}{2}}y^{-\frac{1}{2}}  \int_0^\infty ( \frac{1}{x} \sum_{n=-\infty}^{\infty} e^{- \pi \frac{n^2}{  x^2} - 2 i \pi  z \beta  \frac{n}{x\sqrt{\pi}} })
 (\sum_{p \in \mathbb{Z}} (-1)^{p} e^{ - \pi p^2 x^2 y^2} -1)   x^{-s}dx$

$  =e^{ i \frac{z^2}{2}}y^{-\frac{1}{2}}  \int_0^\infty (x \sum_{n=-\infty}^{\infty} e^{- \pi n^2  x^2 - 2 i \pi  z \beta x \sqrt{\pi} n })  (\sum_{p \in \mathbb{Z}} (-1)^{p} e^{ - \pi p^2 \frac{y^2}{ x^2}} -1)   x^{s-1}dx$

\vspace {0.2 cm}

\textbf{For Operator $\mathbf{H}_{22}$:}
\vspace {0.1 cm}

For this operator a symetric interesting eigenfunction is given by taking $g(t)= \sum_{p \neq 0}^{}  e^{-\frac{p^2}{t^2}} sgn(p)  | p | ^{\frac{1}{2}}$ (and change of variable: $y$ changed by $py$)

$C_s(z,u)= e^{ i  \frac{z^2}{2}} e^{ i  \frac{u^2}{2}} \int_0^\infty \int_0^\infty \sum_{n \neq 0}^{}  (-1)^{n-1}  e^{  (\alpha z-i \sqrt{\pi} n x)^2 } \sum_{p \neq 0}^{} (-1)^{p-1} e^{  (\alpha u-i \sqrt{\pi}p y)^2 } e^{-\frac{k}{x^2y^2}} y^{-\frac{1}{2}}x ^{-s} dx dy $ 

This function is a good candidate as for $s$ zero of Zeta we have the two following relations (which is the $g(t)$ constant term cancellation "effect" under the integrand which we already explained):

$e^{ i  \frac{z^2}{2}} e^{ i  \frac{u^2}{2}}\int_0^\infty  \int_0^\infty \sum_{n \neq 0}^{}  (-1)^{n-1}  e^{  (\alpha z-i \sqrt{\pi} n x)^2 } \sum_{p \neq 0}^{}  (-1)^{p-1} e^{  (\alpha u-i \sqrt{\pi}p y)^2 } y^{-\frac{1}{2}}x ^{-s} dx dy =0 $ 

And the following relation (due to the fact that for Eta function we have a Poisson summation formula):

$C_s(z,u) =e^{ i  \frac{z^2}{2}} e^{ i  \frac{u^2}{2}} \int_0^\infty  \int_0^\infty \sum_{n \neq 0}^{} (-1)^{n-1}  e^{  (\alpha z-i \sqrt{\pi} n x)^2 } \sum_{p \neq 0}^{}(-1)^{p-1}  e^{  (\alpha u-i \sqrt{\pi}p y)^2 } e^{-\frac{k}{x^2y^2}} y^{-\frac{1}{2}}x ^{-s} dx dy$

$ = e^{- i  \frac{z^2}{2}} e^{- i  \frac{u^2}{2}} \int_0^\infty  \int_0^\infty \sum_{n \neq 0}^{}  e^{  (i \alpha z-i \sqrt{\pi} \frac{(2n+1)}{2} x)^2 } \sum_{p \neq 0}^{} e^{  (i \alpha u-i \sqrt{\pi}\frac{(2p+1)}{2} y)^2 } e^{- k x^2y^2} y^{-\frac{1}{2}} x ^{s-1} dx dy $

\vspace {0.2 cm}

Note that if the $g(t)$ functions exists and then that for each zero of Zeta function an eigenfunction for one of the described two dimensional operator, we would have the Hilbert-polya conjecture verified but not in the form of its original idea:
the Hermitian operator alone proposed here would explain only the location of the zeros on the line but not their distribution on this line.

By the way, if $g(t)$ function can be found, then we see that zeros of Zeta and other L-functions with "good" properties (for example Dirichlet L-functions) have their zeros on the "same" line. 

\section{Conclusion} 

\vspace {0.3 cm}

We shown in this article how the L-functions (and especially the Zeta-function) are naturally linked to the operator $\mathbf{H}_{11}= i(x \frac{\partial}{\partial x} - y \frac{\partial}{\partial y})$ which can take the forms $\mathbf{H}_{12}$ or $\mathbf{H}_{22}$ using the transformation  $T(f)(\alpha x)$ on $x$ and/or $y$.

We shown that the interest of such operators is due to the behavior of some of their eigenfunctions which have specific properties when the eigenvalue is the imaginary part of zero of Zeta function (A part of the eigenfunction is canceled).

This cancellation effect is obvious when taking the following eigenfunction of $\mathbf{H}_{11}$ which can be considered as the "initial" Hermitian operator (having $\lambda$ as eigenvalue):

$ x^{-s} y^{-\frac{1}{2}}  \sum_{n=1}^{\infty} \frac{(-1)^{n-1}} {n^{s-1}} g( \frac{x y}{ n}) $

the constant term of the development of $g(t)$ is obviously canceled in this expression if and only if $s$ is zero of the Zeta function (as Eta function as same non trivial zeros as Zeta). 
But we see, on this example, that this specificity will never be enough to have a square module integrable eigenfunction to conclude that the imaginary parts of Zeta function are real. 

To have module integrable functions we proposed to transform this initial operator $\mathbf{H}_{11}$ with the transform $T$ in order to obtain new Hermitian operators having new eigenfunctions for which this "cancellation" effect could be "efficient enough" to put the eigenfunction in $L^2$ .

Note that the operator proposed (taking different forms) is the same as the one proposed in \cite{SO98} (under form proposed p11).

The natural proposal of the operator $x \frac{\partial}{\partial x} - \frac{1}{2}$ is a good idea as this one dimensional operator is naturally associated to Zeta Function, but its eigenfunctions $x^s$ do not present any obvious specificity for $s$ zero of Zeta function and therefore finding a transform of this operator or a condition which would lead to square integrable eigenfunctions is difficult.
The idea to have an initial two dimensional operator allows to have initial eigenfunctions with a clear specificity appearing for $s$ zero of Zeta function.

\vspace {0.3 cm}
We were not able to check if proposed operators could have or not a real spectrum in the proposed Hilbert spaces and if a $g(t)$ functions allowing to have one of the proposed eigenfunction module integrable could exist or not; 
but we consider that the idea to transform the initial operator $\mathbf{H}_{11}$ into a new operator by a transform which will: 

- keep the Hermitian property of the operator 

- transform the initial "cancellation effect" into property for eigenfunctions to be square module integrable.

is an interesting approach to the Riemann Hypothesis.

\vspace {0.5 cm}
Thanks to send your remarks to: $bertrand\_barrau@hotmail.com$

\section{APPENDIX} 

\subsection{Hilbert-Polya conjecture: reminder}   

We remind here that the Hilbert-Polya conjecture is based on the classical following result.

 Noting $s= \frac{1}{2} + i \lambda$, if we have simultaneously:

\begin{itemize}

\item  A Hilbert space $\mathcal{H}$ with a scalar product $ <g , h> = \int g \; \overline{ h}  $

\item A Hermitian operator $H$ defined on the space $\mathcal{H}$

\item A eigenfunction $f_s$ of $H$ such that when $s$ is zero of Zeta function then $f_s$ is in the hilbert space $\mathcal{H}$ and $H f_s = \lambda f_s$. 

\end{itemize}

Then we can use the fact that $H$ is hermitian with the following relation:

$ \lambda <f_s, f_s> = <H f_s,  f_s> = <f_s, H f_s > = \overline{ \lambda } <f_s, f_s> $ 

showing that we have: $\lambda = \overline{\lambda }$ (ie. $\lambda$ real and therefore $s$ is located on the critical line: $Re(s)= \frac{1}{2}$).

The principal difficulty is to find a Hermitian operator "linked" to the Zeta-function.
The link exists if the eigenfunctions $f_s$ of $H$ have something "special" when $s$ is zero of Zeta-function. If this specificity makes the $f_s$ to be in the Hilbert space defined above (so in $L^2$) when $s$ is zero of Zeta-function, then Riemann Hypothesis will be proved.

\subsection{Notations} 

In all the article $s$ designates a complex number in the critical strip (i.e. by definiton in this article: $0 < Re(s) < 1$), when mentionned this $s$ will be considered as a zero of the Zeta function. 

For commodity the variable $s$ will be sometimes replaced by the complex variable $\lambda$ defined by: $s=\textstyle\frac {1}{2}+i \lambda$

When we refer to $\zeta (s)$ zero without precision, we always refer to non trivial zeros (i.e. Zeros having their real part between 0 and 1). 

We note: $\alpha=e^{-i \frac{\pi}{4}}$

When speaking about $L^2$ we speak about module square integrable functions on a two dimensional domain, generally on $(\mathbb{R}^+ \times \mathbb{R})$ or on  $(\mathbb{R} \times \mathbb{R})$.

We consider $L(s)=\displaystyle\sum_{n=1}^{\infty} \dfrac{a_n}{n^s}$ converging in the critical strip and we pose $a_{-n}=a_n$

In this article we consider the following classical Hilbert spaces with their scalar products:
\begin{itemize}

\item  $ <g , h> = \int_{- \infty}^{\infty} \int_{0}^{\infty} g(z,y) \; \overline{ h(z,y)}  dy dz$   on $\mathcal{H}_0 =\{ h(z,y) \in \mathbb{C};  h \in L^2 (\mathbb{R} \times \mathbb{R}^+) ; \forall z, h(z,0)=0\}$

\item  $ <g , h> = \int_{- \infty}^{\infty} \int_{-\infty}^{\infty} g(z,y) \; \overline{ h(z,y)}  dy dz$  on $\mathcal{H}_1 =\{ h(z,y) \in \mathbb{C};  h \in L^2 (\mathbb{R} \times \mathbb{R}) \}$ 

\end{itemize}

Each time we refer to a Hermitian operator we implicitly refer to an operator on $\mathcal{H}_0$ or $\mathcal{H}_1$.

We use classical notations for the Zeta-function $\zeta(s)$, Eta function $\eta(s)$ and the Xi-function $\xi(s)$ : 

\begin{itemize}
\item $\zeta(s)= \displaystyle\sum_{n=1}^{\infty} \dfrac{1}{n^s}$~;
\item $\eta(s)= \displaystyle\sum_{n=1}^{\infty} \dfrac{(-1)^{n-1}}{n^s}$~;
\item $\xi(s)=(s-1) \pi^{-\frac{s}{2}} \Gamma(\frac{s}{2}) \zeta(s)$~;
\end{itemize}

\bibliographystyle{plain}
\bibliography{Biblio1}

\begin{thebibliography}{10}

\bibitem{BK98}
M.V. Berry and J.P. Keating.
\newblock H=xp and the riemann zeros.
\newblock {\em In Supersymmetry and Trace Formulae: Chaos and Disorders}, 1998.

\bibitem{BK99}
M.V. Berry and J.P. Keating.
\newblock The riemann zeros and eigenvalue asymptotics.
\newblock {\em SIAM Review}, 1998.

\bibitem{DB97}
Par~Kurlberg by~Daniel~Bump, Kwok-Kwong~Choi and Jeffrey Vaaler.
\newblock A local riemann hypothesis, i.
\newblock April 1998.

\bibitem{HME74}
Harold~M. Edwards.
\newblock {\em Riemann's zeta Function}.
\newblock Dover Publication - Academic Press, 1974.

\bibitem{YH97}
Yves Hellegouarch.
\newblock {\em Invitation aux Mathematiques de Fermat-Wiles}.
\newblock Dunod edition, 1997.

\bibitem{MAIS72}
M.Abramowitz and I.A. Stegun.
\newblock {\em Handbook Of Mathematical Functions}.
\newblock Dover Publication, 1972.

\bibitem{MO97}
Y.~Motohashi.
\newblock {\em Spectral theory of the Riemann Zeta-Function}.
\newblock Cambridge University Press, 1997.

\bibitem{SO98}
Susumu Okubo.
\newblock Lorentz-invariant hamiltonian and riemann hypothesis.
\newblock {\em J.Phys.A31:1049-1057}, 1998.

\bibitem{BP01}
J.~Pitman P.~Biane and M.~Yor.
\newblock {\em Probability laws related to the Jacobi theta and Riemann zeta
  functions, and Brownian excursions}.
\newblock Bull. Amer. Math. Soc. 38 (2001), no. 4, p. 435 465., 2001.

\bibitem{PAT88}
S.J. Patterson.
\newblock {\em An introduction to the theory of the Riemann Zeta-Function}.
\newblock Cambridge University Press, 1988.

\bibitem{BK95}
A.~Khare R.K.~Bhaduri and J.~Law.
\newblock Phase of the riemann zeta-function and the inverted harmonic
  oscillator.
\newblock {\em Phys. Rev. E 52, 486 (1995), chao-dyn/9406006;}, 1995.

\bibitem{HR03}
Haret~C. Rosu.
\newblock Quantum hamiltonians and prime numbers.
\newblock {\em Mod. Phys. Lett. A 18 (2003) 1205-1213}, May 2003.

\bibitem{GS08}
G.~Sierra and P.~Townsend.
\newblock Landau levels and riemann zeros.
\newblock {\em [arXiv : 0805.4079 math-ph]}, Sep 2008.

\bibitem{TIT86}
E.C. Titchmarsh.
\newblock {\em The theory of the Riemann Zeta-Function}.
\newblock Oxford science Publication, 1986.

\end{thebibliography}

\end{document}